
\magnification\magstephalf
\baselineskip14pt

\input picmac
\input epsf

\def\figfile{ptpats}
\def\fig#1/{\epsfbox{\figfile.#1}}

\def\meno{\medskip\noindent}
\def\pprod{\mathbin{\underline{\mkern-2mu\times\mkern-2mu}}}
\def\display#1:#2:#3\par{\par\hangindent #1 \noindent
	\hbox to #1{\hfill #2 \hskip .1em}\ignorespaces#3\par}
\def\proof{\medskip\noindent{\bf Proof. \ }}
\def\pfbox
  {\hbox{\hskip 3pt\lower2pt\vbox{\hrule
  \hbox to 7pt{\vrule height 7pt\hfill\vrule}
  \hrule}}\hskip3pt}
\def\ddddots{\mathinner{\mskip1mu\raise7pt\vbox{\kern7pt\hbox{.}}\mskip2mu
    \raise4pt\hbox{.}\mskip2mu\raise1pt\hbox{.}\mskip1mu}}
\def\jn{\mathbin{\cdot\mkern-8mu-\mkern-8mu\cdot}}
\def\bib[#1] {\par\noindent\hangindent 40pt\hbox to 20pt{[#1]\hfil}}

\def\Cvet{1}
\def\Frob{2}
\def\GMi{3}
\def\GMii{4}
\def\GMiii{5}
\def\GVL{6}
\def\MaMi{7}

\centerline{\bf Partitioned Tensor Products and Their Spectra}
\centerline{Donald E. Knuth}
\centerline{Computer Science Department,
Stanford University,
Stanford, CA 94305}

\bigskip
{\narrower\smallskip\noindent
{\bf Abstract.}\quad
A pleasant family of graphs defined by Godsil and McKay is shown to have easily
computed eigenvalues in many cases.
\smallskip}

\meno
Let $G$ and $H$ be directed graphs on the respective vertices~$U$ and~$V\!$,
 and
suppose that the vertex sets have each been partitioned into disjoint subsets
$U=U_0\cup U_1$ and $V=V_0\cup V_1$. The {\it partitioned tensor product\/} 
$G\pprod H$ of $G$ and~$H$ with respect to this partitioning is defined as
follows:

\display 20pt:a):
Each vertex of $U_0$ is replaced by a copy of $H\vert V_0$, the subgraph of~$H$
induced by~$V_0$;

\display 20pt:b):
Each vertex of $U_1$ is replaced by a copy of $H\vert V_1$;

\display 20pt:c):
Each arc of $G$ that runs from $U_0$ to $U_1$ is replaced by a copy of the arcs
of~$H$ that run from~$V_0$ to~$V_1$;

\display 20pt:d):
Each arc of $G$ that runs from $U_1$ to $U_0$ is replaced by a copy of the arcs
of~$H$ that run from~$V_1$ to~$V_0$.

\midinsert
$$
\vcenter{\fig1/} \hskip.75in \vcenter{\fig11/}
$$
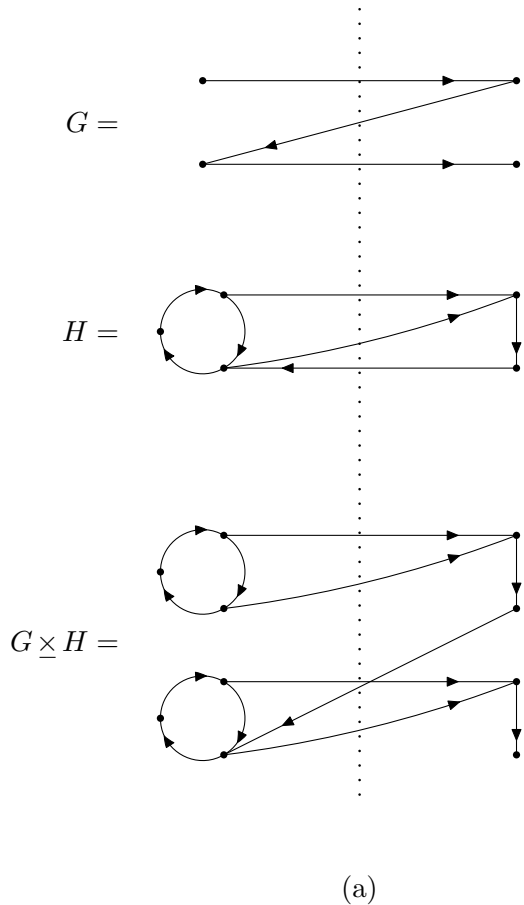
\centerline{Figure 1. Partitioned tensor products, directed and undirected.}
\endinsert

\meno
For example, Figure 1 shows two partitioned tensor products. The example in
Figure~1b is undirected; this is the special case of a directed graph where
each undirected edge corresponds to a pair of arcs in opposite directions.
Arcs of~$G$ that stay within~$U_0$ or~$U_1$ do not contribute to $G\pprod H$,
so we may assume that no such arcs exist (i.e., that $G$ is bipartite).

Figure~2 shows what happens if we interchange
the roles of~$U_0$ and~$U_1$ in~$G$ but leave
everything else intact. (Equivalently, we could interchange the roles of~$V_0$
and~$V_1$.) These graphs, which may be denoted $G^R\pprod H$ to distinguish
them from the graphs  $G\pprod H$  in Figure~1, might look quite different from
their left-right duals,
yet it turns out that the characteristic polynomials of $G\pprod
H$ and $G^R\pprod H$ are strongly related.

\midinsert
$$
\vcenter{\fig2/} \hskip.75in \vcenter{\fig21/}
$$
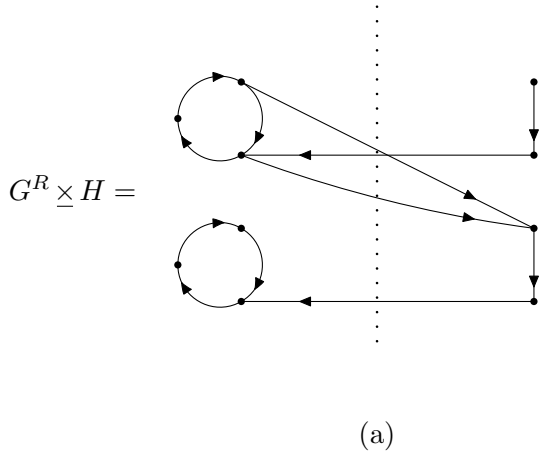
\centerline{Figure 2. Dual products after right-left reflection of $G$.}
\endinsert

Let $E_{ij}$ be the arcs from $U_i$ to $U_j$ in~$G$, and $F_{ij}$ the arcs
from~$V_i$ to $V_j$ in~$H$; multiple arcs are allowed, so $E_{ij}$ and~$F_{ij}$
are multisets. It follows that $G\pprod H$ has $\vert U_0\vert\,\vert V_0\vert
+\vert U_1\vert\,\vert V_1\vert$ vertices and $\vert U_0\vert\,\vert
F_{00}\vert +\vert U_1\vert\,\vert F_{11}\vert +\vert E_{01}\vert\,\vert
F_{01}\vert +\vert E_{10}\vert\,\vert F_{10}\vert$ arcs. Similarly,
$G^R\pprod H$ has $\vert U_1\vert\,\vert V_0\vert+\vert U_0\vert\,\vert
V_1\vert$ vertices and $\vert U_1\vert\,\vert F_{00}\vert+\vert U_0\vert\,
\vert F_{11}\vert +\vert E_{10}\vert\,\vert F_{01}\vert +\vert E_{01}\vert\,
\vert F_{10}\vert$ arcs.

The definition of partitioned tensor product is due to Godsil and McKay
[\GMi],
who proved the remarkable fact that
$$p(G\pprod H)\,p(H\vert V_0)^{\vert U_1\vert-\vert U_0\vert}
=p(G^T\pprod H)\,p(H\vert V_1)^{\vert U_1\vert-\vert U_0\vert}\,,$$
where $p$ denotes the characteristic polynomial of a graph. They also observed
[\GMii] 
that Figures~1b and~2b represent the smallest pair of connected undirected
graphs having the same spectrum (the same~$p$). The purpose of the present note
is to refine their results by showing how to calculate $p(G\pprod H)$
explicitly in terms of~$G$ and~$H$.

We can use the symbols $G$ and $H$ to stand for the adjacency matrices as
well as for the graphs themselves. Thus we have
$$G=\pmatrix{G_{00}&G_{01}\cr G_{10}&G_{11}\cr}\qquad\hbox{and}\qquad
H=\pmatrix{H_{00}&H_{01}\cr H_{10}&H_{11}\cr}$$
in partitioned form, where $G_{ij}$ and $H_{ij}$ denote the respective
adjacency matrices corresponding to the arcs~$E_{ij}$ and~$F_{ij}$. (These
submatrices are not necessarily square; $G_{ij}$~has size $\vert U_i\vert\times
\vert U_j\vert$ and~$H_{ij}$ has size $\vert V_i\vert\times \vert V_j\vert$.)
It follows by definition that
$$G\pprod H=\pmatrix{I_{\vert U_0\vert}\otimes H_{00}&G_{01}\otimes H_{01}\cr
\noalign{\smallskip}
G_{10}\otimes H_{10}&I_{\vert U_1\vert}\otimes H_{11}\cr}$$
where $\otimes$ denotes the Kronecker product or tensor product 
[\MaMi, page~8]
and $I_k$ denotes an identity matrix of size $k\times k$.

Let $H\uparrow\sigma$ denote the graph obtained from $H$ by $\sigma$-fold
repetition of each arc that joins~$V_0$ to~$V_1$. In matrix form
$$H\uparrow\sigma =\pmatrix{\phantom{\sigma} H_{00}&\sigma H_{01}\cr
\sigma H_{10}&\phantom{\sigma} H_{11}\cr}\,.$$
This definition applies to the adjacency matrix
 when $\sigma$ is any complex number, but of course
$H\uparrow\sigma$ is difficult to ``draw'' unless $\sigma$ is a nonnegative
integer. We will show that the characteristic polynomial of $G\pprod H$ factors
into characteristic polynomials of graphs $H\uparrow\sigma$, times a power of
the characteristic polynomials of~$H_{00}$ or~$H_{11}$. The proof is simplest
when $G$ is undirected.

\proclaim
Theorem 1. Let $G$ be an undirected graph, and let $(\sigma_1,\ldots,\sigma_l)$
be the singular values of $G_{01}=G_{10}^T$, where $l=\min(\vert U_0\vert,\vert
U_1\vert)$. Then
$$p(G\pprod H)=\left\{
\vcenter{\halign{$#$\hfil\quad&#\hfil\cr
\left(\prod_{j=1}^l\,p(H\uparrow\sigma_j)\right)
p(H_{00})^{\vert U_0\vert-\vert U_1\vert}\,,&if\/ $\vert U_0\vert\geq
\vert U_1\vert$;\cr
\noalign{\smallskip}
\left(\prod_{j=1}^l\,p(H\uparrow\sigma_j)\right)
p(H_{11})^{\vert U_1\vert-\vert U_0\vert}\,,&if\/ $\vert U_1\vert\geq
\vert U_0\vert$.\cr}}\right.$$

\proof
Any real $m\times n$ matrix $A$ has a singular value decomposition
$$A=QSR^T$$
where $Q$ is an $m\times m$ orthogonal matrix, $R$ is an $n\times n$ orthogonal
matrix, and $S$ is an $m\times n$ matrix with $S_{jj}=\sigma_j\geq 0$ for
$1\leq j\leq \min(m,n)$ and $S_{ij}=0$ for $i\neq j$
[\GVL, page 16].
The numbers $\sigma_1,\ldots,\sigma_{\min(m,n)}$ are called the singular values
of~$A$.

Let $m=\vert U_0\vert$ and $n=\vert U_1\vert$, and suppose that $QSR^T$ is the
singular value decomposition of~$G_{01}$. Then $(\sigma_1,\ldots,\sigma_l)$ are
the nonnegative eigenvalues of the bipartite graph~$G$, and we have
$$\pmatrix{Q^T\otimes I_{\vert V_0\vert}&O\cr 
\noalign{\smallskip}
O&R^T\otimes I_{\vert V_1\vert}\cr}
\,G\pprod H\pmatrix{Q\otimes I_{\vert V_0\vert}&O\cr
\noalign{\smallskip}
 O&R\otimes I_{\vert V_1\vert}\cr}
=\pmatrix{I_{\vert U_0\vert}\otimes H_{00}&S\otimes H_{01}\cr
\noalign{\smallskip}
S^T\otimes H_{10}&I_{\vert U_1\vert}\otimes H_{11}\cr}$$
because $G_{10}=RS^TQ^T$. Row and column permutations of this matrix transform
it into the block diagonal form
$$\pmatrix{H\uparrow\sigma_1&\phantom{0}&\phantom{0}&\phantom{0}\cr
&\ddddots\cr
&&H\uparrow\sigma_l\cr
&&&D\cr}\,,$$
where $D$ consists of $m-n$ copies of $H_{00}$ if $m\geq n$, or $n-m$ copies
of~$H_{11}$ if $n\geq m$. \ \pfbox

A similar result holds when $G$ is directed, but we cannot use the singular
value decomposition because the eigenvalues of~$G$ might not be real and the
elementary divisors of $\lambda I-G$ might not be linear. The following lemma
can be used in place of the singular value decomposition in such cases.

\proclaim
Lemma. Let $A$ and $B$ be arbitrary matrices of complex numbers, where 
$A$ is $m\times n$ and $B$ is $n\times m$. Then we can write
$$A=QSR^{-1}\,,\qquad B=RTQ^{-1}\,,$$
where $Q$ is a nonsingular $m\times m$ matrix, $R$ is a nonsingular $n\times n$
matrix, $S$~is an $m\times n$ matrix, $T$~is an $n\times m$ matrix, and the
matrices $(S,T)$ are triangular with consistent diagonals:
$$\eqalign{S_{ij}&=T_{ij}=0\qquad\hbox{for $i>j$};\cr
S_{jj}&=T_{jj}\quad{\rm or}\quad S_{jj}T_{jj}=0\,,\qquad\hbox{for }
1\leq j\leq \min(m,n)\,.\cr}$$

\proof
We may assume that $m\leq n$. If $AB$ has a nonzero eigenvalue $\lambda$, let
$\sigma$ be any square root of~$\lambda$ and let $x$ be a nonzero $m$-vector
such that $ABx=\sigma^2x$. Then the $n$-vector $y=Bx/\sigma$ is nonzero, and we
have
$$Ay=\sigma x\,,\qquad Bx=\sigma y\,.$$
On the other hand, if all eigenvalues of $AB$ are zero, let $x$ be a nonzero
vector such that $ABx=0$. Then if $Bx\neq 0$, let $y=Bx$. If $Bx=0$, let $y$ be
any nonzero vector such that $Ay=0$; this is possible unless all $n$~columns
of~$A$ are linearly independent, in which case we must have $m=n$ and we can
find $y$ such that $Ay=x$. In all cases we have therefore demonstrated the
existence of nonzero vectors~$x$ and~$y$ such that
$$Ay=\sigma x\,,\qquad Bx=\tau y\,,\qquad \sigma=\tau\quad\hbox{or}\quad
\sigma\tau=0\,.$$

Let $X$ be a nonsingular $m\times m$ matrix whose first column is $x$, and let
$Y$ be a nonsingular $n\times n$ matrix whose first column is~$y$. Then
$$X^{-1}AY=\pmatrix{\sigma&a\cr 0&A_1\cr}\,,\qquad
Y^{-1}BX=\pmatrix{\tau&b\cr 0&B_1\cr}$$
where $A_1$ is $(m-1)\times (n-1)$ and $B_1$ is $(n-1)\times (m-1)$. If $m=1$,
let $Q=X$, $R=Y$, $S=(\sigma a)$, and $T={\tau\choose 0}$. Otherwise we have
$A_1=Q_1S_1R_1^{-1}$ and $B_1=R_1T_1Q_1^{-1}$ by induction, and we can let
$$Q=X\pmatrix{1&0\cr 0&Q_1\cr}\,,\quad R=Y\pmatrix{1&0\cr 0&R_1\cr}\,,\quad
S=\pmatrix{\sigma&aR_1\cr 0&S_1\cr}\,,\quad T=\pmatrix{\tau&BQ_1\cr
0&T_1\cr}\,.$$
All conditions are now fulfilled. \ \pfbox

\proclaim
Theorem 2. Let $G$ be an arbitrary graph, and let $(\sigma_1,\ldots,\sigma_l)$
be such that $\sigma_j=S_{jj}=T_{jj}$ or $\sigma_j=0=S_{jj}T_{jj}$ when
$G_{01}=QSR^{-1}$ and $G_{10}=RTQ^{-1}$ as in the lemma, where $l=\min
(\vert U_0\vert\,,\vert U_1\vert)$. Then $p(G\pprod H)$ satisfies the
identities of Theorem~1.

\proof
Proceeding as in the proof of Theorem 1, we have
$$\pmatrix{Q^{-1}\otimes I_{\vert V_0\vert}&O\cr 
\noalign{\smallskip}
O&R^{-1}\otimes I_{\vert V_1\vert}\cr}
\,G\pprod H\pmatrix{Q\otimes I_{\vert V_0\vert}&O\cr
\noalign{\smallskip}
 O&R\otimes I_{\vert V_1\vert}\cr}
=\pmatrix{I_{\vert U_0\vert}\otimes H_{00}&S\otimes H_{01}\cr
\noalign{\smallskip}
T\otimes H_{10}&I_{\vert U_1\vert}\otimes H_{11}\cr}\,.$$
This time a row and column permutation converts the right-hand matrix to a
block {\it triangular\/} form, with zeroes below the diagonal blocks. Each
block on the diagonal is either $H\uparrow\sigma_j$ or~$H_{00}$
or~$H_{11}$, or of the form
$$\pmatrix{H_{00}&\sigma H_{01}\cr \tau H{10}&H_{11}\cr}\,,\qquad
\sigma\tau =0\,.$$
In the latter case the characteristic polynomial is clearly
$p(H_{00})p(H_{11})=p(H\uparrow 0)$, so the remainder of the proof of Theorem~1
carries over in general. \ \pfbox

\medskip
The proof of the lemma shows that the numbers $\sigma_1^2,\ldots,\sigma^2_p$
are the characteristic roots of $G_{01}G_{10}$, when $\vert U_0\vert\leq\vert 
U_1\vert$, otherwise they are
the characteristic roots of $G_{10}G_{01}$. Either square
root of~$\sigma_j^2$ can be chosen, since the matrix~~$H\uparrow\sigma$ is
similar to~$H\uparrow(-\sigma)$.

We have now reduced the problem of computing $p(G\pprod H)$ to the problem of
computing the characteristic polynomial of the graphs~$H\uparrow\sigma$. The
latter is easy when $\sigma=0$, and some graphs~$G$ have only a few nonzero
singular values. For example, if $G$ is the complete bipartite graph having
parts~$U_0$ and~$U_1$ of sizes~$m$ and~$n$, all singular values vanish except
for $\sigma=\sqrt{mn}$.

If $H$ is small, and if only a few nonzero $\sigma$ need to be considered, the
computation of $p(H\uparrow \sigma)$ can be carried out directly. 
For example, it turns out that
$$\pmatrix{\phantom{-}\lambda&-1&-\sigma&\phantom{-}0&\phantom{-}0\cr
-1&\phantom{-}\lambda&\phantom{-}0&\phantom{-}0&-\sigma\cr
-\sigma&\phantom{-}0&\phantom{-}\lambda&-1&\phantom{-}0\cr
\phantom{-}0&\phantom{-}0&-1&\phantom{-}\lambda&-1\cr
\phantom{-}0&-\sigma&\phantom{-}0&-1&\phantom{-}\lambda\cr}
=(\lambda^2+\lambda-\sigma^2)\bigl(\lambda^3-\lambda^2-
(2+\sigma^2)\lambda+2\bigr)\,;$$
so we can compute the spectrum of $G\pprod H$ by solving a few quadratic and
cubic equations, when $H$ is this particular 5-vertex graph (a~partitioned
5-cycle). 
But it is interesting  to look for large families of graphs for which simple
formulas yield $p(H\uparrow\sigma)$ as a function of~$\sigma$.

One such family consists of graphs that have only one edge crossing the
partition. 
Let $H_{00}$ and~$H_{11}$ be graphs on~$V_0$ and~$V_1$, and form the graph
$H=H_{00}\jn H_{11}$ by adding a single edge between designated vertices
$x_0\in V_0$ and $x_1\in V_1$. Then a glance at the adjacency matrix of~$H$
shows that 
$$p(H\uparrow\sigma)=p(H_{00})p(H_{11})-\sigma^2p(H_{00}\vert V_0\backslash
x_0)p(H_{11}\vert V_1\backslash x_1)\,.$$
(The special case $\sigma=1$ of this formula is Theorem 4.2(ii) of
[\GMiii].)

Another case where $p(H\uparrow\sigma)$ has a simple form arises when the
matrices
$$H_0=\pmatrix{H_{00}&0\cr 0&H_{11}\cr}\qquad{\rm and}\qquad
H_1=\pmatrix{0&H_{01}\cr H_{10}&0\cr}$$
commute with each other. Then it is well known
[\Frob]
that the eigenvalues of $H_0+\sigma H_1$ are $\lambda_j+\sigma\mu_j$, for some
ordering of the eigenvalues~$\lambda_j$ of~$H_0$ and~$\mu_j$ of~$H_1$. Let us
say that $(V_0,V_1)$ is a {\it compatible partition\/} of~$H$ if
$H_0H_1=H_1H_0$, i.e., if
$$H_{00}H_{01}=H_{01}H_{11}\qquad{\rm and}\qquad H_{11}H_{10}=H_{10}H_{00}\,.$$
When $H$ is undirected, so that $H_{00}=H_{00}^T$ and $H_{11}=H_{11}^T$ and
$H_{10}=H_{01}^T$, the compatibility condition boils down to the single
relation 
$$H_{00}H_{01}=H_{01}H_{11}\,.\eqno(\ast)$$

Let $m=\vert V_0\vert$ and $n=\vert V_1\vert$, so that $H_{00}$ is $m\times m$,
$H_{01}$~is $m\times n$, and $H_{11}$ is $n\times n$.
One obvious way to satisfy $(\ast)$ is to let $H_{00}$ and $H_{11}$ both be
zero, so that $H$ is bipartite as well as~$G$. Then $H\uparrow\sigma$ is
simply~$\sigma H$, the $\sigma$-fold repetition of the arcs of~$H$, and its
eigenvalues are just those of~$H$ multiplied by~$\sigma$. For example, if $G$
is the $M$-cube $P_2^M$ and~$H$ is a path~$P_N$ on $N$~points, and if $U_0$
consists of the vertices of even parity in~$G$ while $V_0$ is one of~$H$'s
bipartite parts, the characteristic polynomial of $G\pprod H$ is
$$\prod_{\scriptstyle 1\leq j\leq M\atop\scriptstyle 1\leq k\leq N}\,
\biggl(\lambda-(2N-4j)\cos\;{k\pi\over N+1}\,\biggr)^{{M\choose j}/2}\,,$$
because of the well-known eigenvalues of~$G$ and~$H$
[\Cvet].
Figure~3 illustrates this construction in the special case $M=N=3$. The
smallest pair of cospectral graphs, 
$\,\vcenter{\hbox{\unitlength=5pt
\beginpicture(2,2)(0,0)
\put(0,0){\disk{.3}}
\put(0,2){\disk{.3}}
\put(2,0){\disk{.3}}
\put(2,2){\disk{.3}}
\put(1,1){\disk{.3}}
\put(0,0){\line(1,1)2}
\put(0,2){\line(1,-1)2}
\endpicture}}\,$
and
$\,\vcenter{\hbox{\unitlength=5pt
\beginpicture(2,2)(0,0)
\put(0,0){\disk{.3}}
\put(0,2){\disk{.3}}
\put(2,0){\disk{.3}}
\put(2,2){\disk{.3}}
\put(1,1){\disk{.3}}
\put(0,0){\line(1,0)2}
\put(0,2){\line(1,0)2}
\put(0,0){\line(0,1)2}
\put(2,0){\line(0,1)2}
\endpicture}}\,$,
is obtained in a similar way by considering the eigenvalues of
$P_3\pprod P_3$ and $P_3^T\pprod P_3$
[\GMii].

\midinsert
$$
\vcenter{\fig3/}\hskip.75in\hbox{Figure 3. $P_2^3\pprod P_3$.}
$$
\endinsert

Another simple way to satisfy the compatibility condition $(\ast)$ with
symmetric matrices~$H_{00}$ and~$H_{11}$ is to let $H_{01}$ consist entirely
of~1s, and to let $H_{00}$ and~$H_{11}$ both be regular graphs of the same
degree~$d$. Then the eigenvalues of~$H_0$ are $(\lambda_1,\ldots,\lambda_m,
\lambda'_1,\ldots,\lambda'_n)$, where $(\lambda_1,\ldots,\lambda_m)$ belong
to~$H_{00}$ and $(\lambda'_2,\ldots,\lambda'_n)$ belong to~$H_{11}$ and
$\lambda_1=\lambda'_1=d$. The eigenvalues of~$H_1$ are
$(\sqrt{mn},-\sqrt{mn},0,\ldots,0)$. We can match the eigenvalues of~$H_0$
properly with those of~$H_1$ by looking at the common eigenvectors
$(1,\ldots,1)^T$ and $(1,\ldots,1,-1,\ldots,-1)^T$ that correspond to~$d$
in~$H_0$ and $\pm\sqrt{mn}$ in~$H_1$; the eigenvalues of $H\uparrow\sigma$ are
therefore
$$(d+\sigma\sqrt{mn},\lambda_2,\ldots,\lambda_m,d-\sigma\sqrt{mn},\lambda'_2,
\ldots,\lambda'_n)\,.$$

Yet another easy way to satisfy $(\ast)$ is to assume that $m=n$ and to let
$H_{00}=H_{11}$ commute with~$H_{01}$. One general construction of this kind
arises when the vertices of~$V_0$ and~$V_1$ are the elements of a group, and
when $H_{00}=H_{11}$ is a Cayley graph on that group. In other words, two
elements~$\alpha$ and~$\beta$ are adjacent in~$H_{00}$ iff $\alpha\beta^{-1}\in
X$, where $X$ is an arbitrary set of group elements closed under inverses. And
we can let $\alpha\in V_0$ be adjacent to $\beta\in V_1$ iff
$\alpha\beta^{-1}\in Y$, where $Y$ is any normal subgroup. Then $H_{00}$
commutes with~$H_{01}$. The effect is to make the cosets of~$Y$ fully
interconnected between~$V_0$ and~$V_1$, while retaining a more interesting
Cayley graph structure inside~$V_0$ and~$V_1$. If $Y$ is the trivial subgroup,
so that $H_{01}$ is simply the identity matrix, our partitioned tensor product
$G\pprod H$ becomes simply the ordinary Cartesian product $G\oplus H=I_{\vert
U\vert}\otimes H+G\otimes I_{\vert V\vert}$. But in many other cases this
construction gives something more general.

A fourth family of compatible partitions is illustrated by the following
graph~$H$ in which $m=6$ and $n=12$:
$$\left(\,\,\vcenter{\halign{%
\hfil#\hfil\quad&\hfil#\hfil\quad&\hfil#\hfil\quad&\hfil#\hfil\quad%
&\hfil#\hfil\quad&\hfil#\hfil\quad\hskip3pt
&\hfil#\hfil\quad&\hfil#\hfil\quad&\hfil#\hfil\quad&\hfil#\hfil\quad%
&\hfil#\hfil\quad&\hfil#\hfil\quad&\hfil#\hfil\quad&\hfil#\hfil\quad%
&\hfil#\hfil\quad&\hfil#\hfil\quad&\hfil#\hfil\quad&\hfil#\hfil\cr
\noalign{\vskip2pt}
0&0&1&1&1&0&1&0&0&0&0&0&0&1&0&0&0&0\cr
0&0&0&1&1&1&0&1&0&0&0&0&0&0&1&0&0&0\cr
1&0&0&0&1&1&0&0&1&0&0&0&0&0&0&1&0&0\cr
1&1&0&0&0&1&0&0&0&1&0&0&0&0&0&0&1&0\cr
1&1&1&0&0&0&0&0&0&0&1&0&0&0&0&0&0&1\cr
0&1&1&1&0&0&0&0&0&0&0&1&1&0&0&0&0&0\cr
\noalign{\vskip3pt}
1&0&0&0&0&0&0&0&1&0&1&0&0&0&0&0&1&0\cr
0&1&0&0&0&0&0&0&0&1&0&1&0&0&0&0&0&1\cr
0&0&1&0&0&0&1&0&0&0&1&0&1&0&0&0&0&0\cr
0&0&0&1&0&0&0&1&0&0&0&1&0&1&0&0&0&0\cr
0&0&0&0&1&0&1&0&1&0&0&0&0&0&1&0&0&0\cr
0&0&0&0&0&1&0&1&0&1&0&0&0&0&0&1&0&0\cr
0&0&0&0&0&1&0&0&1&0&0&0&0&0&1&0&1&0\cr
1&0&0&0&0&0&0&0&0&1&0&0&0&0&0&1&0&1\cr
0&1&0&0&0&0&0&0&0&0&1&0&1&0&0&0&1&0\cr
0&0&1&0&0&0&0&0&0&0&0&1&0&1&0&0&0&1\cr
0&0&0&1&0&0&1&0&0&0&0&0&1&0&1&0&0&0\cr
0&0&0&0&1&0&0&1&0&0&0&0&0&1&0&1&0&0\cr
\noalign{\vskip2pt}
}}
\,\,\right)$$
In general, let $C_{2k}$ be the matrix of a cyclic permutation on
$2k$~elements, and let $m=2k$, $n=4k$. Then we obtain a compatible partition if
$$H_{00}=\bigl(C_{2k}^j+C_{2k}^k+C_{2k}^{-j}\bigr)\,,\quad
H_{01}=(I_{2k}\;\;C_{2k})\,,\quad
H_{11}=\pmatrix{C_{2k}^j+C_{2k}^{-j}&C_{2k}^{k+1}\cr
\noalign{\smallskip}
C_{2k}^{k-1}&C_{2k}^j+C_{2k}^{-j}}\,.$$
The $18\times 18$ example matrix is the special case $j=2$, $k=3$. The 
eigenvalues of $H\uparrow\sigma$ in general are
$$\omega^{jl}+\omega^{-jl}+1\,,\qquad
\omega^{jl}+\omega^{-jl}-1+\sqrt{2}\,\sigma\,,\qquad
\omega^{jl}+\omega^{-jl}-1-\sqrt{2}\,\sigma$$
for $0\leq l<2k$, where $\omega=e^{\pi i/k}$.

Compatible partitionings of digraphs are not difficult to construct. But it
would be interesting to find further examples of undirected graphs, without
multiple edges, that have a compatible partition.

\bigskip
\centerline{\bf References}

\bigskip
\bib
[\Cvet]
Drago{\v s} M. Cvetkovi\'c, Michael Doob, and Horst Sachs, {\sl Spectra of
Graphs\/} (New York: Academic Press, 1980).

\bib
[\Frob]
G. Frobenius, ``\"Uber vertauschbare Matrizen,'' {\sl Sitzungsberichte der
K\"oniglich Preu{\ss}ischen Akademie der Wissenschaften zu Berlin\/} (1896),
601--614. Reprinted in his {\sl Gesammelte Abhandlungen\/ \bf 2} (Berlin:
Springer, 1968), 705--718.

\bib
[\GMi]
C. Godsil and B. McKay, ``Products of graphs and their spectra,'' in
{\sl Combinatorial Mathematics IV}, edited by A. Dold and B. Eckmann,
{\sl Lecture Notes in Mathematics\/ \bf560} (1975), 61--72.

\bib
[\GMii]
C. Godsil and B. McKay, ``Some computational results on the spectra of 
graphs,'' in
{\sl Combinatorial Mathematics IV}, edited by A. Dold and B. Eckmann,
{\sl Lecture Notes in Mathematics\/ \bf560} (1975), 73--82.

\bib
[\GMiii]
C. D.  Godsil and B. D. McKay, ``Constructing cospectral graphs,''
{\sl {\AE}quationes Mathematic{\ae}\/ \bf 25} (1982), 257--268.

\bib
[\GVL]
Gene H. Golub and Charles F. Van Loan, {\sl Matrix Computations\/}
(Baltimore: Johns Hopkins University Press, 1983).

\bib
[\MaMi]
Marvin Marcus and Henrik Minc, {\sl A Survey of Matrix Theory and Matrix
Inequalities\/} (Boston: Allyn and Bacon, 1964).

\bye